\documentclass[12pt]{amsart}

\usepackage{amsmath, amsthm, amssymb, amsfonts, verbatim}

\def\0{{\bf 0}}

\def\D{{\mathcal D}}

\def\R{{\mathbb R}}
\def\Z{{\mathbb Z}}

\theoremstyle{plain}

\newtheorem{thm}{Theorem}

\newtheorem{prop}[thm]{Proposition}
\newtheorem{lem}[thm]{Lemma}

\newtheorem{theorem}{Theorem}
\newtheorem{lemma}[thm]{Lemma}

\theoremstyle{definition}
\newtheorem*{definition}{Definition}

\newtheorem*{theorem*}{Theorem}

\theoremstyle{remark}
\newtheorem{remark}[equation]{Remark}

\begin{document}

\raggedbottom

\numberwithin{equation}{section}

%
%
\newcommand{\MarginNote}[1]{
    \marginpar{
        \begin{flushleft}
            \footnotesize #1
        \end{flushleft}
        }
    }
%
%
\newcommand{\NoteToSelf}[1]{
    }

%
%
\newcommand{\Obsolete}[1]{
    }

\newcommand{\Detail}[1]{
    \MarginNote{Detail}
    \skipline
    \hspace{+0.25in}\fbox{\parbox{4.25in}{\small #1}}
    \skipline
    }

\newcommand{\Todo}[1]{
    \skipline \noindent \textbf{TODO}:
    #1
    \skipline
    }

\newcommand{\Comment}[1] {
    \skipline
    \hspace{+0.25in}\fbox{\parbox{4.25in}{\small \textbf{Comment}: #1}}
    \skipline
    }

%
%

\newcommand{\IntTR}
    {\int_{t_0}^{t_1} \int_{\R^d}}

\newcommand{\IntAll}
    {\int_{-\iny}^\iny}

\newcommand{\Schwartz}
    {\ensuremath \Cal{S}}

\newcommand{\SchwartzR}
    {\ensuremath \Schwartz (\R)}

\newcommand{\SchwartzRd}
    {\ensuremath \Schwartz (\R^d)}

\newcommand{\SchwartzDual}
    {\ensuremath \Cal{S}'}

\newcommand{\SchwartzRDual}
    {\ensuremath \Schwartz' (\R)}

\newcommand{\SchwartzRdDual}
    {\ensuremath \Schwartz' (\R^d)}

\newcommand{\HSNorm}[1]
    {\norm{#1}_{H^s(\R^2)}}

\newcommand{\HSNormA}[2]
    {\norm{#1}_{H^{#2}(\R^2)}}

\newcommand{\Holder}
    {H\"{o}lder }

\newcommand{\Holders}
    {H\"{o}lder's }

\newcommand{\Holderian}
    {H\"{o}lderian }

\newcommand{\HolderRNorm}[1]
    {\widetilde{\Vert}{#1}\Vert_r}

\newcommand{\LInfNorm}[1]
    {\norm{#1}_{L^\iny(\Omega)}}

\newcommand{\SmallLInfNorm}[1]
    {\smallnorm{#1}_{L^\iny}}

\newcommand{\LOneNorm}[1]
    {\norm{#1}_{L^1}}

\newcommand{\SmallLOneNorm}[1]
    {\smallnorm{#1}_{L^1}}

\newcommand{\LTwoNorm}[1]
    {\norm{#1}_{L^2(\Omega)}}

\newcommand{\SmallLTwoNorm}[1]
    {\smallnorm{#1}_{L^2}}

\newcommand{\LpNorm}[2]
    {\norm{#1}_{L^{#2}}}

\newcommand{\SmallLpNorm}[2]
    {\smallnorm{#1}_{L^{#2}}}

\newcommand{\lOneNorm}[1]
    {\norm{#1}_{l^1}}

\newcommand{\lTwoNorm}[1]
    {\norm{#1}_{l^2}}

\newcommand{\MsrNorm}[1]
    {\norm{#1}_{\Cal{M}}}

\newcommand{\FTF}
    {\Cal{F}}

\newcommand{\FTR}
    {\Cal{F}^{-1}}

\newcommand{\InvLaplacian}
    {\ensuremath{\widetilde{\Delta}^{-1}}}

\newcommand{\EqDef}
    {\hspace{0.2em}={\hspace{-1.2em}\raisebox{1.2ex}{\scriptsize def}}\hspace{0.2em}}

%
%

%
%

\title
    [Vanishing viscosity for nondecaying vorticity]
    {A finite time result for vanishing viscosity in the plane with nondecaying vorticity}

\author{Elaine Cozzi}
\address{Department of Mathematical Sciences, Carnegie Mellon University}
\curraddr{}
\email{ecozzi@andrew.cmu.edu}

\subjclass{Primary 76D05, 76C99} 
\date{} 


\keywords{Fluid mechanics, Inviscid limit}

\begin{abstract}
Assuming that initial velocity has finite energy and initial vorticity is bounded in the plane, we show that the unique solutions of the Navier-Stokes equations converge to the unique solution of the Euler equations in the $L^{\infty}$-norm uniformly over finite time as viscosity approaches zero.  We also establish a rate of convergence.
\end{abstract}

\maketitle

\section{Introduction}\label{Introduction}
\noindent We consider the Navier-Stokes equations modeling incompressible viscous fluid flow, given by
\begin{align*}
    \begin{matrix}
        (NS) & \left\{
            \begin{matrix}
                \partial_t v_{\nu} + v_{\nu} \cdot \nabla v_{\nu} - \nu \Delta
                 v_{\nu} = - \nabla p_{\nu} \\
                \text{div } v_{\nu} = 0 \\
                v_{\nu}|_{t = 0} = v_{\nu}^0,
            \end{matrix}
            \right.
    \end{matrix}
\end{align*}
and the Euler equations modeling incompressible non-viscous fluid flow, given by
\begin{align*}
    \begin{matrix}
        (E) & \left\{
            \begin{matrix}
                \partial_t v + v \cdot \nabla v = - \nabla p\\
                \text{div } v = 0 \\
                v|_{t = 0} = v^0.
            \end{matrix}
            \right.
    \end{matrix}
\end{align*}
In this paper, we study the vanishing viscosity limit.  The question of vanishing viscosity addresses whether or not a solution $v_{\nu}$ of ($NS$) converges in some norm to a solution $v$ of ($E$) with the same initial data as viscosity tends to $0$.  

The vanishing viscosity problem is closely tied to uniqueness of solutions to the Euler equations, because the methods used to prove uniqueness can often be applied to show vanishing viscosity.  One of the most important uniqueness results in the plane is due to Yudovich.  He establishes in \cite{Yudovich} the uniqueness of a solution $(v,p)$ to ($E$) in the space $C(\R^+; L^2(\R^2))\times L^{\infty}_{loc}(\R^+;L^2(\R^2))$ when $v^0$ belongs to $L^2(\R^2)$ and $\omega^0$ belongs to $L^p(\R^2)\cap L^{\infty}(\R^2)$ for some $p<\infty$.  For this uniqueness class, Chemin proves in \cite{Chemin} that the vanishing viscosity limit holds in the $L^p$-norm uniformly over finite time, and he establishes a rate of convergence.  (In fact, the author only considers the case $p=2$; however, the proof of the result can easily be generalized to any $p<\infty$.)  

In this paper, we consider the case where initial velocity belongs to $L^2(\R^2)$, while initial vorticity is bounded and does not necessarily belong to $L^p(\R^2)$ for any $p<\infty$.  
The existence and uniqueness of solutions to ($E$) with nondecaying vorticity {\em and} nondecaying velocity was proved by Serfati in \cite{Serfati}.  Specifically, the author proves the following theorem.
\begin{theorem}\label{Serfati}
Let $v^0$ and $\omega^0$ belong to $L^{\infty}(\R^2)$, and let $c\in\R$.  For every $T>0$ there exists a unique solution $(v,p)$ to ($E$) in the space $L^{\infty}([0,T]; L^{\infty}(\R^2))\times L^{\infty}([0,T]; C(\R^2))$ with $\omega\in L^{\infty}([0,T]; L^{\infty}(\R^2))$, $p(0)=c$, and with $\frac {p(t,x)}{|x|}\rightarrow 0$ as $|x|\rightarrow \infty$.  
\end{theorem}
The author also shows that the velocity satisfies the estimate
\begin{equation}\label{linfv}
||v(t)||_{L^{\infty}}\leq C||v^0||_{L^{\infty}}e^{C_1||\omega^0||_{L^{\infty}}t}.
\end{equation}

For the main result of this paper, we assume that initial velocity belongs to $L^2(\R^2)$ rather than $L^{\infty}(\R^2)$; however, the boundedness of the initial vorticity combined with finite energy of the initial velocity imply that the initial velocity is bounded (see Lemma \ref{C1} and Remark \ref{C1revised}).  Therefore our solutions to ($E$) satisfy the initial conditions given in Theorem \ref{Serfati}.  We can conclude from Theorem \ref{Serfati} that the velocity remains bounded as time evolves.  Moreover, in Section \ref{A Few Definitions and Technical Lemmas} we show that energy is conserved when $v^0$ is in $L^2(\R^2)$ and $\omega^0$ is in $L^{\infty}(\R^2)$, so that the solution $v$ to ($E$) belongs to $L^{\infty}_{loc}(\R^+; L^2(\R^2)\cap L^{\infty}(\R^2))$.
  
We prove that if $v^0$ belongs to $L^2(\R^2)$ and $\omega^0$ belongs to $L^{\infty}(\R^2)$, then the unique solutions of the Navier-Stokes equations converge to the unique solution of the Euler equations in the $L^{\infty}$-norm uniformly over finite time as viscosity approaches zero (see Theorem \ref{main}).  This result improves upon a previous result in \cite{C} in the sense that here we show convergence for finite time, whereas in \cite{C} we show convergence for short time only.  In \cite{C}, however, we do not assume initial velocity belongs to $L^{2}(\R^2)$; rather we make the weaker assumption that initial velocity is in $L^{\infty}(\R^2)$.  Despite the different assumptions placed on initial velocity, the strategy for proving this result will follow the strategy of \cite{C}.  

In \cite{C} we consider low, middle, and high frequencies of the difference between the solutions to ($NS$) and ($E$) separately.  For both low and high frequencies we are able to show convergence in the $L^{\infty}$-norm on any finite time interval.  For low frequencies, we utilize the structure of mild solutions to ($NS$) as well as the structure of Serfati solutions to ($E$).  For high frequencies, we apply Bernstein's Lemma and Calderon-Zygmund theory.  For the middle frequencies, in \cite{C} we are only able to show that convergence holds for short time.  To do this we localize the frequencies of the vorticity formulations of ($NS$) and ($E$), and we consider the difference of the two resulting equations.  We are then able to reduce the problem to showing that there exists $T_0>0$ such that the vanishing viscosity limit holds in the Besov space ${\dot{B}}^0_{\infty,\infty}(\R^2)$ uniformly on $[0,T]$.  

We apply a similar approach here; however, in this case we use different methods along with the additional assumption that $v^0\in L^2(\R^2)$ to show that the vanishing viscosity limit holds in the ${\dot{B}}^0_{\infty,\infty}$-norm uniformly over any finite time interval.  Combining this new estimate with the estimates from \cite{C} for the low and high frequencies, we are able to extend our previous short time result to a finite time result.

We begin the proof of the main result as in \cite{C}.  Specifically, we define $S_n v=\chi(2^{-n}D)v$, where $\chi$ is identically 1 on the unit ball and 0 outside of the ball of radius two, and we write  
\begin{equation}\label{3terms1}
\begin{split}
&||v_{\nu}-v||_{L^{\infty}([0,T];L^{\infty}(\R^2))} \leq ||S_{-n}(v_{\nu}-v)||_{L^{\infty}([0,T];L^{\infty}(\R^2))}\\
&\qquad+ ||(S_n-S_{-n})(v_{\nu}-v)||_{L^{\infty}([0,T];L^{\infty}(\R^2))}\\
&\qquad+ ||(Id-S_{n})(v_{\nu}-v)||_{L^{\infty}([0,T];L^{\infty}(\R^2))}.
\end{split}
\end{equation}
Letting $n$ be a function of $\nu$ so that as $n$ approaches infinity, $\nu$ tends to $0$, we show that each of the three terms in ($\ref{3terms1}$) converges to $0$ as $n$ approaches infinity.  The estimates of $||S_{-n}(v_{\nu}-v)||_{L^{\infty}([0,T];L^{\infty}(\R^2))}$ and $||(Id-S_{n})(v_{\nu}-v)||_{L^{\infty}([0,T];L^{\infty}(\R^2))}$ both follow from a straightforward application of Bernstein's Lemma.  The first estimate uses the membership of $v_{\nu}$ and $v$ to $L^{\infty}([0,T];L^2(\R^2))$, and the second relies on the membership of $\omega_{\nu}$ and $\omega$ to $L^{\infty}([0,T]; L^{\infty}(\R^2))$.

To estimate the $L^{\infty}$-norm of $(S_n-S_{-n})(v_{\nu}-v)$, we initially follow the approach given in \cite{C}.  We observe that this term has Fourier support in an annulus with inner and outer radius of order $2^{-n}$ and $2^n$, respectively.  Breaking this support into $2n$ smaller annuli, we use the definition of the Besov space ${\dot{B}}^{0}_{\infty,\infty}$ (see Section \ref{A Few Definitions and Technical Lemmas}) to conclude that 
\begin{equation}\label{intro}
||(S_n-S_{-n})(v_{\nu}-v)(t)||_{L^{\infty}}\leq Cn||(v_{\nu}-v)(t)||_{{\dot{B}}^{0}_{\infty,\infty}}
\end{equation} 
for each $t>0$.  We proceed to show that $v_{\nu}$ converges to $v$ in the ${\dot{B}}^{0}_{\infty,\infty}$-norm on any finite time interval as viscosity approaches $0$ (see Lemma \ref{key}).  It is in this proof that we deviate from the methods of \cite{C}.  Specifically, we replace an argument using Gronwall's Lemma with an argument relying on Osgood's Lemma.  We closely follow the proof in \cite{Chemin} which demonstrates that the vanishing viscosity limit holds in the energy norm uniformly over finite time when initial velocity belongs to $L^2(\R^2)$ and initial vorticity belongs to $L^2(\R^2)\cap L^{\infty}(\R^2)$.  

\Obsolete{By reducing the problem of vanishing viscosity in the $L^{\infty}$-norm to that in the ${\dot{B}}^{0}_{\infty,\infty}$-norm, we eliminate some of the general difficulties one encounters when proving estimates for the fluid equations in $L^{\infty}(\R^n)$.  For example, Calderon-Zygmund operators are not bounded on $L^{\infty}(\R^n)$ but are bounded on ${\dot{B}}^{s}_{\infty,\infty}(\R^n)$ for all $s\in\R$ (see \cite{Torres}).  We can then use the relation $\nabla v=\nabla{\nabla}^{\perp}{\Delta}^{-1}\omega(v)$, the boundedness of Calderon-Zygmund operators on ${\dot{B}}^{s}_{\infty,\infty}$, and Lemma \ref{homoprop} below to see that the ${\dot{B}}^{0}_{\infty,\infty}$-norm of $v$ and the ${\dot{B}}^{-1}_{\infty,\infty}$-norm of $\omega(v)$ are equivalent.  This equivalence allows us to prove the main result by using the vorticity formulations of ($NS$) and ($E$) to show that the vorticities $\omega(v_{\nu})$ converge to $\omega(v)$ in the ${\dot{B}}^{-1}_{\infty,\infty}$-norm as viscosity tends to $0$. }

The paper is organized as follows.  In Section \ref{A Few Definitions and Technical Lemmas}, we introduce the Littlewood-Paley theory, Bony's paraproduct decomposition, Besov spaces, and some useful technical lemmas.  In Section \ref{Properties of solutions to the euler equations with nondecaying vorticity}, we discuss properties of solutions to ($E$) with nondecaying vorticity.  In Section \ref{Statement and Proof of the Main Result} and Section \ref{Proof of Lemma}, we state and prove the main result; we devote Section \ref{Proof of Lemma} entirely to showing that the vanishing viscosity limit holds in the ${\dot{B}}^0_{\infty,\infty}$-norm.
\section{A Few Definitions and Technical Lemmas}\label{A Few Definitions and Technical Lemmas}
We first define the Littlewood-Paley operators (see \cite{Triebel} for further discussion of these operators).  We let $\varphi \in S(\R^d)$ satisfy supp $\varphi\subset \{\xi\in \R^d: \frac{3}{4} \leq|\xi |\leq \frac{8}{3} \}$, and for every $j\in\Z$ we let $\varphi_j(\xi)=\varphi(2^{-j}\xi)$ (so
${\check{\varphi}}_j(x)=2^{jd}\check{\varphi}(2^jx))$.  Observe that, if $|j-j'|\geq 2$, then supp ${\varphi_j}$ $\cap$ supp
${\varphi_{j'}} = \emptyset$.  We define ${\psi}_n \in S(\R^d)$ by the equality
\begin{equation*}
{\psi}_n (\xi) = 1 - \sum_{j\geq n} \varphi_j(\xi) 
\end{equation*}
for all $\xi\in\R^d$, and for $f\in S'(\R^d)$ we define the operator $S_nf$ by  
\begin{equation*}
S_n f = {\check{\psi}}_n \ast f.
\end{equation*}
In the following sections we will make frequent use of both the homogeneous and the inhomogeneous Littlewood-Paley operators.  For $f\in S'(\R^d)$ and $j\in\Z$, we define the homogeneous Littlewood-Paley operators ${\dot{\Delta}}_j$ by
\begin{equation*}
{\dot{\Delta}}_j f= {\check{\varphi}}_j \ast f,
\end{equation*}
and we define the inhomogeneous Littlewood-Paley operators by
\begin{align*}
    \Delta_j f = \left\{
        \begin{array}{ll}
            0,
                & j < -1, \\
            \check{\psi_0} \ast f,
                & j = -1, \\
            {\check{\varphi}}_j \ast f,
                & j > -1.
        \end{array}
        \right.
\end{align*}  
We remark that the operators $\Delta_j$ and $\dot{\Delta}_j$ coincide when $j\geq 0$, but differ when $j\leq -1$.

In the proof of the main theorem we use the paraproduct decomposition introduced by J.-M. Bony in \cite{Bony}.
We recall the definition of the paraproduct and remainder used in this decomposition. 
\begin{definition}\label{para}
Define the paraproduct of two functions $f$ and $g$ by
\begin{equation*}
T_fg = \sum_{\stackrel{i,j}{i\leq j-2}} \Delta_i f\Delta_j g = \sum_{j=1}^{\infty} S_{j-1}f\Delta_j g.
\end{equation*}
We use $R(f,g)$ to denote the remainder.  $R(f,g)$ is given by the following bilinear operator:  
\begin{equation*}
R(f,g)=\sum_{\stackrel{i,j}{|i-j|\leq 1}} \Delta_if\Delta_jg.
\end{equation*}
\end{definition}
\noindent Bony's decomposition then gives
\begin{equation*}
fg=T_fg + T_gf + R(f,g).
\end{equation*} 
We now define the homogeneous Besov spaces.
\begin{definition}\label{besovhomo}
Let $s\in\R$, $(p,q)\in[1,\infty]\times[1,\infty)$.  We define the homogeneous Besov space ${\dot{B}}^s_{p,q}(\R^d)$ to be the space of tempered distributions $f$ on $\R^d$ such that\\
\begin{equation*}
||f||_{{\dot{B}}^s_{p,q}}:={\left(\sum_{j\in\Z}2^{jqs}{||{\dot{\Delta}}_jf||}^q_{L^p}\right)}^{\frac{1}{q}} < \infty.
\end{equation*} 
When $q=\infty$, write
\begin{equation*}
||f||_{{\dot{B}}^s_{p,\infty}}:={\sup_{j\in\Z}{2^{js}{||{\dot{\Delta}}_jf||}_{L^p}}}.
\end{equation*} 
\end{definition}
\Obsolete{\begin{definition}\label{besov}
Let $s\in\R$, $(p,q)\in[1,\infty]\times[1,\infty)$.  We define the inhomogeneous Besov space $B^s_{p,q}(\R^d)$ to be the space of tempered distributions $f$ on $\R^d$ such that
\begin{equation*}
||f||_{B^s_{p,q}}:=||S_0f||_{L^p}+{\left(\sum_{j=0}^{\infty}2^{jqs}{||{\dot{\Delta}}_jf||}^q_{L^p}\right)}^{\frac{1}{q}} < \infty.
\end{equation*} 
When $q=\infty$, write
\begin{equation*}
||f||_{{B}^s_{p,\infty}}:=||S_0f||_{L^p}+{\sup_{j\geq 0}{2^{js}{||{\dot{\Delta}}_jf||}_{L^p}}}.
\end{equation*} 
\end{definition}}
We also define the inhomogeneous Zygmund spaces.
\begin{definition}
The inhomogeneous Zygmund space $C^s_{\ast}(\R^d)$ is the set of all tempered distributions $f$ on $\R^d$ such that 
\begin{equation*}
||f||_{{C}^s_{\ast}}:=\sup_{j\geq -1}2^{js}||{\Delta}_j f||_{L^{\infty}} < \infty.
\end{equation*} 
\end{definition}
It is well-known that $C^s_{\ast}(\R^d)$ coincides with the classical Holder space $C^s(\R^d)$ when $s$ is not an integer and $s>0$ (see for example \cite{Chemin1}).

We will make frequent use of Bernstein's Lemma.  We refer the reader to \cite{Chemin1}, chapter 2, for a proof of the lemma.
\begin{lem}\label{bernstein}
(Bernstein's Lemma) Let $r_1$ and $r_2$ satisfy $0<r_1<r_2<\infty$, and let $p$ and $q$ satisfy $1\leq p \leq q \leq \infty$. There exists a positive constant $C$ such that for every integer $k$ , if $u$ belongs to $L^p(\R^d)$, and supp $\hat{u}\subset B(0,r_1\lambda)$, then 
\begin{equation}\label{bern1}
\sup_{|\alpha|=k} ||\partial^{\alpha}u||_{L^q} \leq C^k{\lambda}^{k+d(\frac{1}{p}-\frac{1}{q})}||u||_{L^p}.
\end{equation}
Furthermore, if supp $\hat{u}\subset C(0, r_1\lambda, r_2\lambda)$, then 
\begin{equation}\label{bern2}
C^{-k}{\lambda}^k||u||_{L^p} \leq \sup_{|\alpha|=k}||\partial^{\alpha}u||_{L^p} \leq C^{k}{\lambda}^k||u||_{L^p}.
\end{equation} 
\end{lem} 
As a result of Bernstein's Lemma, we have the following lemma regarding the homogeneous Besov spaces.
\begin{lem}\label{homoprop}
Let $u\in S'(\R^d)$, $s\in\R$, and $p,q\in[1,\infty]$.  For any $k\in\Z$  there exists a constant $C_k$ such that whenever $|\alpha|=k$,
\begin{equation*}
{C}^{-k}||\partial^{\alpha}u||_{{\dot{B}}^s_{p,q}} \leq ||u||_{{\dot{B}}^{s+k}_{p,q}} \leq C^k||\partial^{\alpha}u||_{{\dot{B}}^s_{p,q}}.
\end{equation*}
\end{lem}
We also make use of the following technical lemma.  We refer the reader to \cite{CK} for a detailed proof.
\begin{lemma}\label{CZhighfreq}
Let $v$ be a divergence-free vector field with vorticity $\omega$, and let $v$ and $\omega$ satisfy the relation $\nabla v=\nabla{\nabla}^{\perp}{\Delta}^{-1}\omega$.  Then there exists an absolute constant $C$ such that for all $j\in\Z$,
    \begin{align*}
        ||{\dot{\Delta}}_j \nabla v||_{L^{\infty}}
            \leq C||{\dot{\Delta}}_j \omega||_{L^{\infty}}.
    \end{align*}
\end{lemma} 
\section{Properties of solutions to the Euler equations with nondecaying vorticity}\label{Properties of solutions to the euler equations with nondecaying vorticity}
In this section, we discuss properties of solutions to ($E$) with nondecaying vorticity.  We begin by proving a result giving Holder regularity of solutions to ($E$) with initial velocity and initial vorticity in $L^{\infty}(\R^2)$.  Specifically, we prove that under these assumptions on the initial data, the corresponding solution to ($NS$) or ($E$) belongs to the Zygmund space $C^{1}_{\ast}(\R^2)$.
\begin{lemma}\label{C1}
Let $v$ be the unique solution to ($E$) given by Theorem \ref{Serfati} with bounded initial vorticity and bounded initial velocity.  Then the following estimate holds:
\begin{equation*}
||v(t)||_{C^{1}_{\ast}} \leq C||v^0||_{L^{\infty}}e^{C_1||\omega^0||_{L^{\infty}}t} + C||\omega^0||_{L^{\infty}}.
\end{equation*}
\end{lemma}
\begin{proof}
Write 
\begin{equation*}
||v(t)||_{C^{1}_{\ast}} \leq C||S_0 v(t)||_{L^{\infty}} + \sup_{j\geq 0}2^j||{\dot{\Delta}}_j v(t)||_{L^{\infty}}.
\end{equation*}
We first use Young's inequality combined with ($\ref{linfv}$) to bound the low frequency term.  For the high frequency terms, we apply Bernstein's Lemma and the uniform estimate on the $L^{\infty}$-norm of the vorticity given in ($\ref{NSvortbound}$) below to bound the supremum by $C||\omega^0||_{L^{\infty}}$.  This completes the proof. 
\end{proof}

We now address the question of energy conservation when $\omega^0$ belongs to $L^{\infty}(\R^2)$ and $v^0$ belongs to $L^2(\R^2)$.  As mentioned in Section \ref{Introduction}, with these assumptions on the initial data, Theorem \ref{Serfati} guarantees that the solution $v$ to ($E$) belongs to $L^{\infty}_{loc}(\R^+; L^{\infty}(\R^2))$.   However, in \cite{Serfati} the author does not assume $v^0\in L^2(\R^2)$ and therefore does not address the issue of energy conservation when vorticity is nondecaying.  For this we prove the following lemma.
\begin{lemma}\label{energy}
Assume $v^0$ belongs to $L^2(\R^2)$ and $\omega^0$ belongs to $L^{\infty}(\R^2)$.  Then the solution to ($E$) given by Theorem \ref{Serfati} satisfies $||v(t)||_{L^2}=||v^0||_{L^2}$ for all $t>0$.
\end{lemma}
\begin{proof}   
In \cite{Serfati}, when the author proves existence of a solution ($v,p$) to ($E$) satisfying the assumptions of Theorem \ref{Serfati}, he assumes $v^0$ and $\omega^0$ belong to $L^{\infty}(\R^2)$, and he uses the sequence $\{S_nv^0\}$ to generate a sequence of smooth solutions $\{(v_n, p_n)\}$ to ($E$) with initial velocity $S_nv^0$.  The author then proceeds to prove that the sequence $\{v_n\}$ is uniformly bounded in $W^{1,p}_{loc}([0,T]\times\R^2)$ for $p\in(2+\epsilon, \infty)$, $\epsilon>0$.  Using a standard diagonalization argument, one can construct a subsequence, which we henceforth denote as $\{v_n\}$, converging to $v$ in $W^{1,q}_{loc}([0,T]\times\R^2)$ for $q<p$.  To complete the proof, the author passes to the limit in ($E$) written in weak form to show that $(v,p)$ solves ($E$).  (We refer the reader to \cite{Serfati} for the details of this argument.)  

To prove Lemma \ref{energy}, we use the same sequence $\{v_n\}$ of smooth solutions as that used in \cite{Serfati}, but we now assume that $v^0$ also belongs to $L^2(\R^2)$.  In this case we have, for each $t>0$ and for each $n$,
\begin{equation*}
||v_n(t)||_{L^2}=||S_nv^0||_{L^2}\leq C||v^0||_{L^2},
\end{equation*}  
where $C$ is an absolute constant.  Therefore, for each fixed $t$, there exists a subsequence $\{v_{n_k}(t)\}$ converging weakly in $L^2(\R^2)$ to $\bar{v}(t)$ with $||\bar{v}(t)||_{L^2}\leq C||v^0||_{L^2}$.  Moreover, by the argument in \cite{Serfati}, $\{v_{n_k}\}$ converges to the solution $v$ of ($E$) in $L^2_{loc}([0,T]\times\R^2)$.  This implies $\{v_{n_k}(t)\}$ converges to $v(t)$ weakly in $L^2(\Omega)$ for each bounded set $\Omega\subset \R^2$ and for almost every $t\in[0,T]$.  Using uniqueness of weak limits we can conclude that $\bar{v}(t)=v(t)$ for almost every $t$ in $[0,T]$; thus the weak solution $v$ generated from initial data $v^0$ belongs to $L^{\infty}_{loc}(\R^+; L^{2}(\R^2)\cap L^{\infty}(\R^2))$.     

To complete the proof of Lemma \ref{energy}, we use a result proved by Duchon and Robert in \cite{DR}.  In \cite{DR}, the authors show that weak solutions $v$ to ($E$) in $L^3_{loc}(\R^+; L^3(\R^2))$ conserve energy if
\begin{equation}\label{energy2}
D_{\epsilon}(v)(t,x)=\frac{1}{4}\int_{\R^2}{\nabla {\phi}^{\epsilon} (y) \cdot \delta v(t) (\delta v(t))^2} dy
\end{equation}
approaches $0$ in $\D'((0,T)\times\R^2)$ as $\epsilon$ approaches $0$, where $\phi$ is a standard mollifier on $\R^2$, ${\phi}^{\epsilon} (y)=\frac{1}{\epsilon^2}\phi(\frac{y}{\epsilon})$, and $\delta v(t)=v(t,x+y)-v(t,x)$.
 

\Obsolete{To prove the lemma, we use an argument identical to that found in \cite{DR}.  In \cite{DR}, the authors show that weak solutions $v$ to ($E$) in $L^3_{loc}(\R^+; L^3(\R^2))$ conserve energy if
\begin{equation}\label{energy2}
D_{\epsilon}(v)(x)=\frac{1}{4}\int{\nabla {\phi}^{\epsilon} (y) \cdot \delta v (\delta v)^2} dy
\end{equation}
approaches $0$ in the sense of distributions as $\epsilon$ approaches $0$, where $\phi$ is a standard mollifier on $\R^2$, ${\phi}^{\epsilon} (y)=\frac{1}{\epsilon^2}\phi(\frac{y}{\epsilon})$, and $\delta v=v(x+y)-v(x)$.  We do not know a priori that $v$ belongs to $L^3_{loc}(\R^+; L^3(\R^2))$.  We do, however, have sufficient regularity of the pair $(v,p)$ to ensure that energy is still conserved for our weak solutions if ($\ref{energy2}$) approaches $0$ as $\epsilon$ approaches $0$.  

To see this, we follow the proof given in \cite{DR}, making modifications where necessary.  As in \cite{DR} we let $u^{\epsilon}=\phi^{\epsilon}\ast u$ and we regularize the Euler equation which gives
\begin{equation}\label{1978}
\partial_t v^{\epsilon} +\partial_i(v_iv)^{\epsilon} +\nabla p^{\epsilon} =0.
\end{equation}
Taking the dot product of ($\ref{1978}$) with $v$ and adding the resulting equation to ($E$) multiplied by $v^{\epsilon}$ yields
\begin{equation}\label{1979}
\partial_t(v\cdot v^{\epsilon}) + \nabla\cdot \left((v\cdot v^{\epsilon})v + v p^{\epsilon} + pv^{\epsilon}\right) + E_{\epsilon} = 0,
\end{equation}           
where 
\begin{equation}\label{1981}
E_{\epsilon}=\partial_i((v_iv_j)^{\epsilon}v_j)-v_iv_j\partial_iv_j^{\epsilon}.
\end{equation}
Since $v$ belongs to the space $L^{\infty}_{loc}(\R^+; C^1_{\ast}(\R^2))$ by Lemma \ref{C1}, it follows that $v\cdot v^{\epsilon}$ converges to $v^2$ in $\D'([0,T]\times \R^2)$ as $\epsilon$ approaches $0$.  Moreover, using the membership of $v$ to $L^{\infty}_{loc}(\R^+; C^1_{\ast}(\R^2))$ and the membership of $p$ to $L^{\infty}_{loc}(\R^+;C(\R^2))$ as given by Theorem \ref{Serfati}, we see that $(v\cdot v^{\epsilon})v + v p^{\epsilon} + pv^{\epsilon}$ converges to $(v^2+2p)v$ in $\D'([0,T]\times\R^2)$.  We can conclude from ($\ref{1979}$) that $E_{\epsilon}$ converges in the sense of distributions to 
\begin{equation*}
-\partial_t(v^2) - \nabla\cdot (v(v^2+2p)).
\end{equation*}  
As in \cite{DR}, it follows by a calculation that
\begin{equation}\label{1980}
\int{\nabla {\phi}^{\epsilon} (y) \cdot \delta v (\delta v)^2} dy = -\partial_i(v_iv_jv_j)^{\epsilon} + 2\partial_i(v_iv_j)^{\epsilon}v_i+\partial_i (v_jv_j)^{\epsilon}v_i-2v_iv_j\partial_iv_j^{\epsilon}.
\end{equation} 
Moreover, by incompressibility of $v$ we can conclude that $\partial_i (v_jv_j)^{\epsilon}v_i=\partial_i ((v_jv_j)^{\epsilon}v_i)$.  Since $v$ belongs to $L^{\infty}_{loc}(\R^+; C^1_{\ast}(\R^2))$, it follows that $v_jv_j$ and $v_iv_jv_j$ are bounded and uniformly continuous.  Therefore, the quantity $-\partial_i(v_iv_jv_j)^{\epsilon} + \partial_i (v_jv_j)^{\epsilon}v_i$ converges to $0$ in $\D'(\R^2)$.  Comparing ($\ref{1980}$) and ($\ref{1981}$), we see that $\int{\nabla {\phi}^{\epsilon} (y) \cdot \delta v (\delta v)^2} dy$ and $2E_{\epsilon}$ have the same limit in the sense of distributions.  That is, $\frac{1}{4}\int{\nabla {\phi}^{\epsilon} (y) \cdot \delta v (\delta v)^2} dy$ converges to 
\begin{equation}\label{localenergy}
-\partial_t(\frac{1}{2}v^2) - \nabla\cdot (v(\frac{1}{2}v^2+p))
\end{equation}
in $\D'(\R^2)$ as $\epsilon$ approaches $0$. }

To see that $\frac{1}{4}\int_{\R^2}{\nabla {\phi}^{\epsilon} (y) \cdot \delta v(t) (\delta v(t))^2} dy$ approaches $0$ for the solutions in Lemma \ref{energy}, we utilize the membership of $v$ to $L^{\infty}_{loc}(\R^+; C^1_{\ast}(\R^2))$.  Taking into account the support of $\phi$, we can write 
\begin{equation*}
\begin{split}
\left|\int_{\R^2}{\nabla {\phi}^{\epsilon} (y) \cdot \delta v(t) (\delta v(t))^2} dy \right|&\leq C\int_{\R^2}{\frac{\epsilon}{|y|^{2}}\left|\frac{1}{\epsilon^2}{(\nabla \phi)}(\frac{y}{\epsilon}) \cdot \delta v(t) (\delta v(t))^2\right|}dy \\
&\leq C\epsilon||v(t)||^3_{C^{\frac{2}{3}}}\leq (Ce^{C_1t})\epsilon,
\end{split}
\end{equation*}
where we used Lemma \ref{C1} to get the last inequality, and $C$ and $C_1$ depend only on the initial data.  This completes the proof of Lemma \ref{energy}. 
\end{proof}
\begin{remark}\label{C1revised}
We remark here that we can revise the proof of Lemma \ref{C1} to show that 
\begin{equation}\label{energyforNS}
||v(t)||_{C^1_{\ast}} \leq C||v^0||_{L^2}+C||\omega^0||_{L^{\infty}}.
\end{equation}
To do this we simply apply Gronwall's Lemma and Young's inequality to the low frequencies, and use energy conservation to get the series of inequalities $||S_0v(t)||_{L^{\infty}} \leq C||v(t)||_{L^{2}} = C||v^0||_{L^{2}}$.  We also point out that an estimate identical to ($\ref{energyforNS}$) holds for $v_{\nu}$ by the argument given in the proof of Lemma \ref{C1} combined with the fact that the energy of Leray solutions to ($NS$) is nonincreasing over time.
\end{remark}

Finally, we will need a uniform bound in time on the $L^{\infty}$-norms of the vorticities corresponding to the solutions of ($NS$) and ($E$).  For fixed $\nu\geq 0$, we have that 
\begin{equation}\label{NSvortbound}
||\omega_{\nu}(t)||_{L^{\infty}} \leq ||\omega_{\nu}^0||_{L^{\infty}}
\end{equation}
for all $t\geq 0$.  One can prove this bound by applying the maximum principle to the vorticity formulations of ($NS$) and ($E$).  We refer the reader to Lemma 3.1 of \cite{ST} for a detailed proof.  

\section{Statement and Proof of the Main Result}\label{Statement and Proof of the Main Result}
\noindent We are now prepared to state the main theorem.
\begin{theorem}\label{main}
Let $v_{\nu}$ be the unique Leray solution to ($NS$) and let $v$ be the unique solution to ($E$), both with initial data $v^0$ in $L^2(\R^2)$ and $\omega^0$ in $L^{\infty}(\R^2)$.  Then there exist constants $C$ and $C_1$, depending only on the initial data, such that the following estimate holds for any fixed $T>0$ and for any $\alpha\in (0,1)$:
\begin{equation}\label{thebigone}
||v_{\nu}-v||_{L^{\infty}([0,T];L^{\infty}(\R^2))} \leq C(T+1)(\sqrt{\nu})^{\alpha e^{-C_1T}}.
\end{equation}
\end{theorem} 
\Obsolete{\begin{remark}
The constant $C_1$ in ($\ref{thebigone}$) is equal to $A||\omega^0||_{L^{\infty}}$, where $A$ is an absolute constant.  Therefore, one can conclude from Theorem \ref{main} that the vanishing viscosity limit holds on a time interval with length inversely proportional to the size of $||\omega^0||_{L^{\infty}}$.  Moreover, the smallness of $\nu$ required to conclude ($\ref{thebigone}$) depends on $T$ (and therefore on $||\omega^0||_{L^{\infty}}$).  Specifically, larger $T$ requires smaller $\nu$ for ($\ref{thebigone}$) to hold (see Remark \ref{killnu}). 
\end{remark}} 
\begin{proof}
To establish the result, we apply a similar strategy to that used to prove Theorem 3 in \cite{C}.  We let $v_{\nu}$ and $v$ be the unique solutions to ($NS$) and ($E$), respectively, satisfying the assumptions of Theorem \ref{main}.  In what follows, we let $v_n=S_nv$ and $\omega_n=S_n\omega(v)$.  We have the following inequality:
\begin{equation}\label{3terms}
\begin{split}
&||v_{\nu}-v||_{L^{\infty}([0,T];L^{\infty}(\R^2))} \leq ||S_{-n}(v_{\nu}-v)||_{L^{\infty}([0,T];L^{\infty}(\R^2))}\\
&\qquad+ ||(Id-S_{-n})(v_{\nu}-v_n)||_{L^{\infty}([0,T];L^{\infty}(\R^2))}\\
&\qquad+ ||(Id-S_{-n})(v_{n}-v)||_{L^{\infty}([0,T];L^{\infty}(\R^2))}.
\end{split}
\end{equation}
We will estimate each of the three terms on the right hand side of the inequality in ($\ref{3terms}$).  We begin with the third term.  We use the definition of $v_n$, Bernstein's Lemma, Lemma \ref{CZhighfreq}, and (\ref{NSvortbound}) to obtain the inequality 
\begin{equation}\label{term3}
||(Id-S_{-n})(v_{n}-v)||_{L^{\infty}([0,T];L^{\infty}(\R^2))} \leq C2^{-n}||\omega^0||_{L^{\infty}}
\end{equation} 
for any $T>0$.  To bound the first term on the right hand side of ($\ref{3terms}$), we use the fact that the Fourier transform of $S_{-n}(v_{\nu}-v)$ is supported on a ball of radius $2^{-n}$, along with Lemma \ref{energy} and Bernstein's Lemma, to conclude that
\begin{equation}\label{term1}
||S_{-n}(v_{\nu}-v)(t)||_{L^{\infty}} \leq C2^{-n}||S_{-n}(v_{\nu}-v)(t)||_{L^2}\leq C2^{-n}||v^0||_{L^2}
\end{equation}
for every $t>0$.  It remains to bound the second term on the right hand side of ($\ref{3terms}$), given by $||(Id-S_{-n})(v_{\nu}-v_n)||_{L^{\infty}([0,T];L^{\infty}(\R^2))}$.  We prove the following proposition.
\begin{prop}\label{technical}
Let $v_{\nu}$ and $v$ be solutions to ($NS$) and ($E$), respectively, satisfying the properties of Theorem \ref{main}.  Then there exist constants $C$ and $C_1$, depending only on the initial data, such that the following estimate holds for any fixed $\alpha\in(0,1)$ and $T>0$, and for sufficiently large $n$:
\begin{equation}\label{zozozo}
||(Id-S_{-n})(v_{\nu}-v_n)||_{L^{\infty}([0,T]; L^{\infty}(\R^2))} \leq C(T+1)2^{-n\alpha e^{-C_1T}}.
\end{equation}
\end{prop}
\begin{remark}
In the proof of Proposition \ref{technical}, we let $\nu=2^{-2n}$.  Therefore, the dependence of the right hand side of ($\ref{zozozo}$) on $\nu$ is hidden in its dependence on $n$. 
\end{remark}
\begin{proof}
We begin with the series of inequalities
\begin{equation}\label{zoey}
\begin{split}
&||(Id-S_{-n})(v_{\nu}-v_n)(t)||_{L^{\infty}} \leq ||(Id-S_{n})(v_{\nu}-v_n)(t)||_{L^{\infty}}\\
&\qquad\qquad +||(S_n-S_{-n})(v_{\nu}-v_n)(t)||_{L^{\infty}}\\
&\qquad\qquad\leq C2^{-n}||\omega^0||_{L^{\infty}} + Cn||(v_{\nu}-v_n)(t)||_{{\dot{B}}^0_{\infty,\infty}}.\\
\end{split}
\end{equation}
The second inequality follows from an application of Bernstein's Lemma, Lemma \ref{CZhighfreq}, and ($\ref{NSvortbound}$).  

To bound $||(v_{\nu}-v_n)(t)||_{{\dot{B}}^0_{\infty,\infty}}$, we need the following lemma, whose proof we postpone until the next section.
\begin{lem}\label{key}
Let $v_{\nu}$ and $v$ be solutions to ($NS$) and ($E$), respectively, satisfying the properties of Theorem \ref{main}.  Then there exist constants $C$ and $C_1$, depending only on the initial data, such that the following estimate holds for any fixed $\alpha\in(0,1)$ and $T>0$:
\begin{equation*}
||v_{\nu}-v_n||_{L^{\infty}([0,T]; {\dot{B}}^0_{\infty,\infty}(\R^2))} \leq C(T+1)2^{-n\alpha e^{-C_1T}}
\end{equation*}
\end{lem}
Since $\alpha\in(0,1)$ is arbitrary, we can write
\begin{equation*}
n||v_{\nu}-v_n||_{L^{\infty}([0,T]; {\dot{B}}^0_{\infty,\infty}(\R^2))} \leq C(T+1)2^{-n\alpha e^{-C_1T}}
\end{equation*}
for sufficiently large $n$.  Combining this estimate with the estimate given in ($\ref{zoey}$) yields Proposition \ref{technical}.  
\end{proof}
To complete the proof of Theorem \ref{main}, we combine ($\ref{3terms}$), ($\ref{term1}$), Proposition \ref{technical}, and ($\ref{term3}$) to get the following estimate for fixed $T>0$ and for sufficiently large $n$:
\begin{equation*}
||(v_{\nu}-v)||_{L^{\infty}([0,T];L^{\infty}(\R^2))} \leq C(T+1)2^{-n\alpha e^{-C_1T}}. 
\end{equation*}
Using the equality $n=-\frac{1}{2}\log_2 {\nu}$, we obtain ($\ref{thebigone}$).  This completes the proof of Theorem \ref{main}.  We devote the next section to the proof of Lemma \ref{key}.  \end{proof}
\section{Proof of Lemma \ref{key}}\label{Proof of Lemma}
We begin with some notation.  For the proof of Lemma \ref{key}, we let ${\bar{\omega}}_n=\omega_{\nu}-\omega_n$ and ${\bar{v}}_n=v_{\nu}-v_n$.

We recall the following estimate, which is proved in Section 5 of \cite{C}:
\Obsolete{localize the frequencies of the vorticity formulations of ($E$) and ($NS$), and we consider the difference of the two resulting equations.  After localizing the frequency of the vorticity formulation of ($E$), we see that ${\dot{\Delta}}_j\omega_n$ satisfies the following equation: 
\begin{equation}\label{1}
\partial_t{\dot{\Delta}}_j\omega_n + v_n\cdot\nabla{\dot{\Delta}}_j\omega_n + [{\dot{\Delta}}_j, v_n\cdot\nabla]\omega_n= \nabla\cdot {\dot{\Delta}}_j{\tau}_n(v,\omega),
\end{equation} 
where  
\begin{equation*}
\tau_n(v,\omega)= r_n(v,\omega) - (v-v_n)(\omega-\omega_n),
\end{equation*}
and
\begin{equation*}
r_n(v,\omega) = \int{\check{\psi}(y)(v(x-2^{-n}y)-v(x))(\omega(x-2^{-n}y)-\omega(x))}dy.
\end{equation*}  
This equation is utilized by Constantin and Wu in \cite{CW} and by Constantin, E, and Titi in a proof of Onsager's conjecture in \cite{CET}.    

If $v_{\nu}$ is a solution to ($NS$), we can localize the frequency of the vorticity formulation of ($NS$) to see that ${\dot{\Delta}}_j\omega_{\nu}$ satisfies
\begin{equation}\label{NS}
\partial_t {\dot{\Delta}}_j\omega_{\nu} + v_{\nu}\cdot\nabla{\dot{\Delta}}_j\omega_{\nu} + [{\dot{\Delta}}_j, v_{\nu}\cdot\nabla]\omega_{\nu} = \nu\Delta{\dot{\Delta}}_j\omega_{\nu}.
\end{equation} 
We subtract ($\ref{1}$) from ($\ref{NS}$).  This yields
\Obsolete{\begin{equation*}
\begin{split}
&\partial_t{\dot{\Delta}}_j{\bar{\omega}}_n + v_n\cdot\nabla{\dot{\Delta}}_j{\bar{\omega}}_n + {\dot{\Delta}}_j({\bar{v}}_n\cdot\nabla {\omega}_{\nu}) -\nu\Delta {\dot{\Delta}}_j {\bar{\omega}}_n -\nu\Delta{{\dot{\Delta}}_j\omega_n} \\
&\qquad + \nabla\cdot {\dot{\Delta}}_j{\tau}_n(v,\omega)+ [{\dot{\Delta}}_j, v_n\cdot\nabla]{\bar{\omega}}_n=0,
\end{split}
\end{equation*}
which gives}
\begin{equation}\label{zozo}
\begin{split}
&\partial_t{\dot{\Delta}}_j{\bar{\omega}}_n + v_n\cdot\nabla{\dot{\Delta}}_j{\bar{\omega}}_n -\nu\Delta {\dot{\Delta}}_j{\bar{\omega}}_n = -{\dot{\Delta}}_j({\bar{v}}_n\cdot\nabla {\omega}_{\nu}) \\
&\qquad + \nu\Delta{{\dot{\Delta}}_j\omega_n} - \nabla\cdot {\dot{\Delta}}_j{\tau}_n(v,\omega) - [{\dot{\Delta}}_j, v_n\cdot\nabla]{\bar{\omega}}_n.
\end{split}
\end{equation}
We observe that $v_n$ is a divergence-free Lipschitz vector field and apply the following lemma, which is proved in \cite{H}. 
\begin{lem}\label{Hmidi}
Let $p\in[1,\infty]$, and let $u$ be a divergence-free vector field belonging to $L^1_{loc}(\R^+; Lip(\R^d))$.  Moreover, assume the function $f$ belongs to $L^1_{loc}(\R^+; L^p(\R^d))$ and the function $a^0$ belongs to $L^p(\R^d)$.  Then any solution $a$ to the problem
\begin{align*}
    \begin{matrix}
         & \left\{
            \begin{matrix}
                \partial_t a + u \cdot \nabla a - \nu \Delta
                 a = f, \\
                a|_{t = 0} = a^0
            \end{matrix}
            \right.
    \end{matrix}
\end{align*}
satisfies the following estimate:
\begin{equation*}
||a(t)||_{L^p}\leq ||a^0||_{L^p} + \int_0^t{||f(s)||_{L^p}}ds.
\end{equation*} 
\end{lem} 
An application of Lemma \ref{Hmidi} to ($\ref{zozo}$) yields
\begin{equation*}
\begin{split}
&||{\dot{\Delta}}_j{\bar{\omega}}_n(t)||_{L^{\infty}} \leq ||{\dot{\Delta}}_j{{\bar{\omega}}_n}^0||_{L^{\infty}} + C\int_0^t  {||(-{\dot{\Delta}}_j({\bar{v}}_n\cdot\nabla {\omega}_{\nu}) + \nu\Delta{{\dot{\Delta}}_j\omega_n} }\\
 &\qquad -{\nabla\cdot {\dot{\Delta}}_j{\tau}_n(v,\omega) - [{\dot{\Delta}}_j, v_n\cdot\nabla]{\bar{\omega}}_n)(s)||_{L^{\infty}}}ds.
\end{split}
\end{equation*}   
Multiplying through by $2^{-j}$ and taking the supremum over $j\in\Z$, we obtain the inequality
}       
\begin{equation}\label{terms2}
\begin{split}
&\partial_t ||{\bar {v}}_n(t)||_{{\dot{B}}^{0}_{\infty,\infty}} \leq C( ||{\bar{v}}_n{\omega}_{\nu}(t)||_{{\dot{B}}^{0}_{\infty,\infty}} + ||\nu\nabla\omega_n(t)||_{{\dot{B}}^{0}_{\infty,\infty}}\\
&\qquad\qquad +  ||{\tau}_n(v,\omega)(t)||_{{\dot{B}}^{0}_{\infty,\infty}} + \sup_{j\in\Z}2^{-j}||[{\dot{\Delta}}_j, v_n\cdot\nabla]{\bar{\omega}}_n(t)||_{L^{\infty}}),
\end{split}
\end{equation} 
where
\begin{equation*}
\tau_n(v,\omega)= r_n(v,\omega) - (v-v_n)(\omega-\omega_n),
\end{equation*}
and
\begin{equation*}
r_n(v,\omega) = \int{\check{\psi}(y)(v(x-2^{-n}y)-v(x))(\omega(x-2^{-n}y)-\omega(x))}dy.
\end{equation*}  
We observe that 
\begin{equation*}
||{\bar {v}}_n(t)||_{{\dot{B}}^{0}_{\infty,\infty}} \leq C||{\bar {v}}_n(t)||_{L^{\infty}}\leq  C(||v^0||_{L^{2}}+||\omega^0||_{L^{\infty}}) 
\end{equation*}
by Remark \ref{C1revised}, and we define $A=1+ C(||v^0||_{L^{2}}+||\omega^0||_{L^{\infty}}) $ with $C$ as in Remark \ref{C1revised}.  If we define $\delta_n(t)$ by 
\begin{equation*}
\delta_n(t)=\frac{||{\bar {v}}_n(t)||_{{\dot{B}}^{0}_{\infty,\infty}} }{A} \leq 1,
\end{equation*}
then, dividing both sides of ($\ref{terms2}$) by $A$, we see that 
\begin{equation}\label{terms3}
\begin{split}
&\partial_t \delta_n(t) \leq \frac{C}{A}(||{\bar{v}}_n{\omega}_{\nu}(t)||_{{\dot{B}}^{0}_{\infty,\infty}} + ||\nu\nabla\omega_n(t)||_{{\dot{B}}^{0}_{\infty,\infty}}\\
&\qquad\qquad +  ||{\tau}_n(v,\omega)(t)||_{{\dot{B}}^{0}_{\infty,\infty}} + \sup_{j\in\Z}2^{-j}||[{\dot{\Delta}}_j, v_n\cdot\nabla]{\bar{\omega}}_n(t)||_{L^{\infty}}).
\end{split}
\end{equation}
To complete the proof of Lemma \ref{key}, we estimate each term on the right hand side of ($\ref{terms3}$).  To bound $||{\bar{v}}_n{\omega}_{\nu}(t)||_{{\dot{B}}^{0}_{\infty,\infty}}$,  we observe that $L^{\infty}$ is continuously embedded in ${\dot{B}}^0_{\infty,\infty}$ and that the $L^{\infty}$-norm of vorticity is uniformly bounded in time by ($\ref{NSvortbound}$) to write 
\begin{equation}\label{badterm}
\begin{split}
&||{\bar{v}}_n{\omega}_{\nu}(t)||_{{\dot{B}}^0_{\infty,\infty}} \leq C||{\bar{v}}_n(t)||_{L^{\infty}}||{\omega}^0||_{L^{\infty}}.
\end{split}
\end{equation}
To estimate the $L^{\infty}$-norm of ${\bar{v}}_n(t)$, we use ($\ref{term1}$), Bernstein's Lemma, Lemma \ref{CZhighfreq}, and ($\ref{NSvortbound}$) to write
\begin{equation}\label{log}
\begin{split}
&||{\bar{v}}_n(t)||_{L^{\infty}} \leq ||S_{-n}(v_{\nu}-v_n)(t)||_{L^{\infty}} + ||(S_n-S_{-n})(v_{\nu}-v_n)(t)||_{L^{\infty}}\\
&+||(Id-S_n)(v_{\nu}-v_n)(t)||_{L^{\infty}} \leq C2^{-n}+C||(S_n-S_{-n})(v_{\nu}-v_n)(t)||_{L^{\infty}}.\\ 
\end{split}
\end{equation}
Combining ($\ref{badterm}$) and ($\ref{log}$) gives 
\begin{equation}\label{Besov}
\begin{split}
&||{\bar{v}}_n {\omega}_{\nu}(t)||_{{\dot{B}}^{0}_{\infty,\infty}} \leq C2^{-n}+C||(S_n-S_{-n})(v_{\nu}-v_n)(t)||_{L^{\infty}}.\\
\end{split}
\end{equation}
To bound $\nu||\nabla{\omega_n}(t)||_{{\dot{B}}^{0}_{\infty,\infty}}$, we again use Bernstein's Lemma, the definition of $\omega_n$, and ($\ref{NSvortbound}$) to conclude that
\begin{equation*}
\begin{split}
&\nu||\nabla{\omega_n}(t)||_{{\dot{B}}^{0}_{\infty,\infty}} 
\leq C\nu 2^{n} ||{\omega}^0||_{L^{\infty}}.
\end{split}
\end{equation*} 
If we let $\nu=2^{-2n}$, we obtain the inequality
\begin{equation}\label{Besov1}
\nu||\nabla{\omega_n}(t)||_{{\dot{B}}^{0}_{\infty,\infty}}\leq C2^{-n} ||{\omega}^0||_{L^{\infty}}.
\end{equation}  
Finally, we estimate $||{\tau}_n(v,\omega)(t)||_{{\dot{B}}^{0}_{\infty,\infty}}$.  We begin by estimating $||(v-v_n)(\omega-\omega_n)(t)||_{{\dot{B}}^{0}_{\infty,\infty}}$.  We again use the embedding $L^{\infty}\hookrightarrow {\dot{B}}^0_{\infty,\infty}$, Bernstein's Lemma, Lemma \ref{CZhighfreq}, and ($\ref{NSvortbound}$) to infer that 
\begin{equation}\label{nonterm1}
\begin{split}
&||(v-v_n)(\omega-\omega_n)(t)||_{{\dot{B}}^{0}_{\infty,\infty}}
\leq ||(v-v_n)(\omega-\omega_n)(t)||_{L^{\infty}}\\
&\qquad\qquad\leq C||\omega^0||_{L^{\infty}}2^{-n}||\omega^0||_{L^{\infty}}.
\end{split}
\end{equation}
In order to bound $||r_n(v,\omega)(t)||_{{\dot{B}}^{0}_{\infty,\infty}}$, we use the membership of $v$ to $L_{loc}^{\infty}(\R^+; C^{\alpha}(\R^2))$ for any $\alpha\in(0,1)$, as given by Lemma \ref{C1}, to write
\begin{equation}\label{useholder}
|v(t,x-2^{-n}y)-v(t,x)| \leq C2^{-n\alpha}|y|^{\alpha}||v(t)||_{C^{\alpha}}.
\end{equation}
Since ${|y|}^{\alpha}|\check{\psi}(y)|$ is integrable, we can apply ($\ref{useholder}$) and Holder's inequality to estimate $r_n(v,\omega)$.  We conclude that
\begin{equation}\label{nonterm2}
||r_n(v,\omega)(t)||_{{\dot{B}}^{0}_{\infty,\infty}} \leq C2^{-n\alpha}||\omega^0||_{L^{\infty}}||v(t)||_{C^{\alpha}}\leq C2^{-n\alpha},
\end{equation}  
where we used Remark \ref{C1revised} to get the last inequality.  Combining ($\ref{nonterm1}$) and ($\ref{nonterm2}$) yields
\begin{equation}\label{Besov2}
||{\tau}_n(v,\omega)(t)||_{{\dot{B}}^{0}_{\infty,\infty}}  \leq C2^{-n\alpha}.
\end{equation}
Combining ($\ref{terms3}$) with ($\ref{Besov}$), ($\ref{Besov1}$), and ($\ref{Besov2}$) gives
\begin{equation}\label{chemin}
\begin{split}
&\partial_t \delta_n(t) \leq \frac{C}{A}(2^{-n\alpha} + ||(S_n-S_{-n})(v_{\nu}-v_n)(t)||_{L^{\infty}}\\
&\qquad\qquad+\sup_{j\in\Z}2^{-j}||[{\dot{\Delta}}_j, v_n\cdot\nabla]{\bar{\omega}}_n(t)||_{L^{\infty}}),
\end{split}
\end{equation}
where $C$ depends only on the initial data.  To estimate $||(S_n-S_{-n})(v_{\nu}-v_n)(t)||_{L^{\infty}}$, we let $p$, $q\in(1,\infty)$ satisfy $\frac{1}{p}+\frac{1}{q}=1$, and we observe that 
\begin{equation}\label{middlefreqest}
\begin{split}
&||(S_n-S_{-n})(v_{\nu}-v_n)(t)||_{L^{\infty}} \leq \sum_{j=-n}^0 || {\dot{\Delta}}_j({\bar{v}}_n)(t)||_{L^{\infty}}^{\frac{1}{p}}|| {\dot{\Delta}}_j({\bar{v}}_n)(t)||_{L^{\infty}}^{\frac{1}{q}}\\
&\qquad \qquad+ \sum_{j=1}^n || {\dot{\Delta}}_j({\bar{v}}_n)(t)||_{L^{\infty}}^{\frac{1}{p}}|| {\dot{\Delta}}_j({\bar{v}}_n)(t)||_{L^{\infty}}^{\frac{1}{q}}\\
&\qquad\leq ||{\bar{v}}_n(t)||_{{\dot{B}}^{0}_{\infty,\infty}}^{\frac{1}{q}}\left(\sum_{j=-n}^0 2^{\frac{j}{p}} || {\dot{\Delta}}_j({\bar{v}}_n)(t)||_{L^2}^{\frac{1}{p}} +  \sum_{j=1}^n 2^{\frac{-j}{p}}||{\bar{v}}_n(t)||_{{\dot{B}}^{1}_{\infty,\infty}}^{\frac{1}{p}}\right)\\
&\leq C||{\bar{v}}_n(t)||_{{\dot{B}}^{0}_{\infty,\infty}}^{\frac{1}{q}}\frac{1}{1-2^{-\frac{1}{p}}} \left(||v^0||_{L^2}^{\frac{1}{p}}+||{\bar{v}}_n(t)||_{{\dot{B}}^{1}_{\infty,\infty}}^{\frac{1}{p}}\right)\leq Cp||{\bar{v}}_n(t)||_{{\dot{B}}^{0}_{\infty,\infty}}^{\frac{1}{q}},
\end{split}
\end{equation}
where we used Bernstein's Lemma to get the second inequality and Lemma \ref{energy} to get the third inequality.  For the last inequality, we bounded $||{\bar{v}}_n(t)||_{{\dot{B}}^{1}_{\infty,\infty}}$ with $C||\omega^0||_{L^{\infty}}$ using Lemma \ref{homoprop}, Lemma \ref{CZhighfreq} and ($\ref{NSvortbound}$), and we bounded $\frac{1}{1-2^{-\frac{1}{p}}}$ with $Cp$, where $C$ is an absolute constant.  Combining ($\ref{middlefreqest}$) and ($\ref{chemin}$) gives
\begin{equation*}
\begin{split}
&\partial_t \delta_n(t) \leq \frac{C}{A}\left(2^{-n\alpha} + p||{\bar{v}}_n(t)||_{{\dot{B}}^{0}_{\infty,\infty}}^{\frac{1}{q}}+\sup_{j\in\Z}2^{-j}||[{\dot{\Delta}}_j, v_n\cdot\nabla]{\bar{\omega}}_n(t)||_{L^{\infty}}\right)\\
&\qquad\leq \frac{C}{A}2^{-n\alpha} + \frac{Cp}{A^{\frac{1}{p}}}\delta_n(t)^{\frac{1}{q}}+\frac{C}{A}\sup_{j\in\Z}2^{-j}||[{\dot{\Delta}}_j, v_n\cdot\nabla]{\bar{\omega}}_n(t)||_{L^{\infty}}.
\end{split}
\end{equation*}
Since $A\geq 1$, we can write
\begin{equation}\label{keyineq}
\begin{split}
&\partial_t \delta_n(t) \leq  C2^{-n\alpha} + Cp\delta_n(t)^{\frac{1}{q}}\\
&\qquad\qquad +\frac{C}{A}\sup_{j\in\Z}2^{-j}||[{\dot{\Delta}}_j, v_n\cdot\nabla]{\bar{\omega}}_n(t)||_{L^{\infty}}.
\end{split}
\end{equation}
It remains to bound the commutator term on the right hand side of ($\ref{keyineq}$).  We prove the following lemma.
\begin{lemma}\label{comm}
Let $v$ be the unique solution to ($E$) satisfying the conditions of Theorem \ref{main}.  Then the following estimate holds for any pair $(p,q)\in(1,\infty)\times(1,\infty)$ satisfying $\frac{1}{p} +\frac{1}{q}=1$:\\
\begin{equation}\label{commest}
\begin{split}
&{\sup_{j\in\Z}2^{-j}||[{\dot{\Delta}}_j, v_n\cdot\nabla]{\bar{\omega}}_n(t)||_{L^{\infty}}}\leq C\left(2^{-n}+p||{\bar{v}}_n(t)||^{\frac{1}{q}}_{{\dot B}^{0}_{\infty,\infty}}\right),
\end{split}
\end{equation}
where $C$ depends only on $||v^0||_{L^{\infty}}$ and $||\omega^0||_{L^{\infty}}$.
\end{lemma}
\begin{proof}
Throughout the proof, to simplify notation we will often omit the time variable $t$.  When $t$ is omitted from a calculation, it is assumed that $t$ is fixed throughout that calculation. 

We first use Bony's paraproduct decomposition to write
\begin{equation}\label{Bony}
\begin{split}
&[{\dot{\Delta}}_j, v_n\cdot\nabla]{\bar{\omega}}_n = \sum_{m=1}^2 [{\dot{\Delta}}_j, T_{v_n^m}\partial_m]{\bar{\omega}}_n\\
&\qquad + [{\dot{\Delta}}_j, T_{\partial_m\cdot}v_n^m]{\bar{\omega}}_n + [{\dot{\Delta}}_j, \partial_mR(v_n^m, \cdot)]{\bar{\omega}}_n.
\end{split}
\end{equation} 
To bound the $L^{\infty}$-norm of the first term on the right hand side of ($\ref{Bony}$), we consider the cases $j<0$ and $j\geq0$ separately.  For $j\geq0$, we use the definition of the paraproduct and properties of the partition of unity to establish the following equality:
\begin{equation}\label{C67}
[{\dot{\Delta}}_j, T_{v_n}\partial_m] = \sum_{j'=\max\{1,j-4\}}^{j+4} [{\dot{\Delta}}_j, S_{j'-1}(v_n)] \Delta_{j'}\partial_m.
\end{equation}
We then express the operator ${\dot{\Delta}}_j$ as a convolution with ${\check{\varphi}}_j$, write out the commutator on the right hand side of ($\ref{C67}$), and change variables.  This yields
\begin{equation}\label{stupidterm1}
\begin{split}
&||[{\dot{\Delta}}_j, T_{v_n}\partial_m]{\bar{\omega}}_n||_{L^{\infty}} \leq \sum_{j'=\max\{1,j-4\}}^{j+4} ||\int {{\check{\varphi}}(y)(S_{j'-1}v_n(x-2^{-j}y)}\\
&\qquad {-S_{j'-1}v_n(x))\Delta_{j'}\partial_m{\bar{\omega}}_n(x-2^{-j}y)}dy||_{L^{\infty}}\\
&\qquad\leq \sum_{j'=\max\{1,j-4\}}^{j+4}2^{j'-j} ||S_{j'-1}\nabla v_n||_{L^{\infty}}|| \Delta_{j'}{\bar{\omega}}_n||_{L^{\infty}}\int {|\check{\varphi}(y)||y|} dy\\
&\qquad\leq C||{\bar{v}}_n||^{\frac{1}{q}}_{{\dot B}^{0}_{\infty,\infty}} \sum_{j'=\max\{1,j-4\}}^{j+4} \left(2^{j'-j}2^{j'}|| \Delta_{j'}{\bar{v}}_n||^{\frac{1}{p}}_{L^{\infty}}\sum_{l=-1}^{j'-2}||\Delta_l\nabla v_n||_{L^{\infty}}\right).
\end{split}
\end{equation}
We bound $||\Delta_l\nabla v_n||_{L^{\infty}}$ by $||v||_{L^{2}}+||\omega||_{L^{\infty}}$ for each $l$ using Remark \ref{C1revised} and we bound $|| \Delta_{j'}{\bar{v}}_n||^{\frac{1}{p}}_{L^{\infty}}$ by $2^{\frac{-j'}{p}}|| \Delta_{j'}{\bar{\omega}}_n||^{\frac{1}{p}}_{L^{\infty}}$ using Bernstein's Lemma.  We can then write
\begin{equation}
\begin{split}
&||[{\dot{\Delta}}_j, T_{v_n}\partial_m]{\bar{\omega}}_n||_{L^{\infty}} \leq ||{\bar{v}}_n||^{\frac{1}{q}}_{{\dot B}^{0}_{\infty,\infty}}(||v||_{L^{2}}+||\omega||_{L^{\infty}})\\
& \qquad\times \sum_{j'=\max\{1,j-4\}}^{j+4}\left( 2^{j'-j} 2^{j'}||\Delta_{j'}{\bar{\omega}}_n||^{\frac{1}{p}}_{L^{\infty}}\sum_{l\leq j'}2^{\frac{l-j'}{p}}2^{\frac{-l}{p}}\right)\\
&\qquad\leq C2^j||{\bar{v}}_n||^{\frac{1}{q}}_{{\dot B}^{0}_{\infty,\infty}}(||v||_{L^{2}}+||\omega||_{L^{\infty}})||{\bar{\omega}}_n||^{\frac{1}{p}}_{L^{\infty}}\frac{1}{1-2^{-\frac{1}{p}}}.
\end{split}
\end{equation}   
We bound $||v||_{L^{2}}$, $||\omega||_{L^{\infty}}$, and $||{\bar{\omega}}_n||^{\frac{1}{p}}_{L^{\infty}}$ using Lemma \ref{energy} and ($\ref{NSvortbound}$), and we bound $\frac{1}{1-2^{-\frac{1}{p}}}$ by $Cp$, where $C$ is an absolute constant.  We then multiply by $2^{-j}$ and take the supremum over $j\geq0$, which gives
\begin{equation*}
\sup_{j\geq0}2^{-j} ||[{\dot{\Delta}}_j, T_{v_n}\partial_m]{\bar{\omega}}_n(t)||_{L^{\infty}} \leq Cp||{\bar{v}}_n(t)||^{\frac{1}{q}}_{{\dot B}^{0}_{\infty,\infty}}.
\end{equation*}
For the case $j<0$, we apply a different strategy.  We first reintroduce the sum over $m$ and expand the commutator. We then use the properties of our partition of unity, the assumption that $j<0$, and the divergence-free assumption on $v$ to establish the following series of inequalities:
\begin{equation}\label{J23}
\begin{split}
&\sum_{m=1}^2 ||[{\dot{\Delta}}_j, T_{v_n}\partial_m] {\bar{\omega}}_n||_{L^{\infty}} \leq \sum_{m=1}^2 \sum_{j'=1}^2 ||\partial_m\dot{\Delta}_j(S_{j'-1}v_n\Delta_{j'}{\bar{\omega}}_n)||_{L^{\infty}}\\
&\qquad\qquad + \sum_{m=1}^2 \sum_{j'=1}^2||S_{j'-1}v_n\Delta_{j'}{\dot{\Delta}}_j\partial_m{\bar{\omega}}_n||_{L^{\infty}} \\
&\qquad\qquad\leq C\sum_{j'=1}^2 2^j ||v||_{L^{\infty}}||\Delta_{j'}{\bar{\omega}}_n||_{L^{\infty}},\\
\end{split}
\end{equation}
where we applied Bernstein's Lemma to get the factor of $2^j$ in the last inequality.  We bound $||\Delta_{j'}{\bar{\omega}}_n||_{L^{\infty}}$ with $2^{j'}||\Delta_{j'}{\bar{v}}_n||_{L^{\infty}}$, again by Bernstein's Lemma, we multiply ($\ref{J23}$) by $2^{-j}$, and we take the supremum over $j<0$.  This gives
\begin{equation*}
\sup_{j<0} 2^{-j}\sum_{m=1}^2 ||[{\dot{\Delta}}_j, T_{v_n}\partial_m]{\bar{\omega}}_n ||_{L^{\infty}} \leq C||v||_{L^{\infty}}||{\bar{v}}_n||_{L^{\infty}}.
\end{equation*}
We now bound $||v||_{L^{\infty}}$ using Remark \ref{C1revised}.  We estimate the $L^{\infty}$-norm of ${\bar{v}}_n$ as in ($\ref{log}$) and we apply ($\ref{middlefreqest}$) to conclude that
\begin{equation*}
\sup_{j<0}\sum_{m=1}^2 2^{-j}||[{\dot{\Delta}}_j, T_{v_n}\partial_m] {\bar{\omega}}_n(t)||_{L^{\infty}} \leq C2^{-n}+Cp||{\bar{v}}_n(t)||_{{\dot{B}}^{0}_{\infty,\infty}}^{\frac{1}{q}}.
\end{equation*}
We now estimate the $L^{\infty}$-norm of $[{\dot{\Delta}}_j, T_{\partial_m\cdot}v_n]{\bar{\omega}}_n$.  We write out the commutator and estimate the $L^{\infty}$-norm of ${\dot{\Delta}}_j(T_{\partial_m {\bar{\omega}}_n}v_n)$ and $T_{\partial_m{\dot{\Delta}}_j {\bar{\omega}}_n}v_n$ separately.  By the definition of the paraproduct and by properties of our partition of unity, we have
\begin{equation}\label{comm2}
\begin{split}
&||T_{\partial_m{\dot{\Delta}}_j {\bar{\omega}}_n}v_n||_{L^{\infty}} = ||\sum_{l\geq 1}S_{l-1}\partial_m{\dot{\Delta}}_j{\bar{\omega}}_n\Delta_l v_n||_{L^{\infty}}\\
&\qquad\leq\sum_{l=\max\{1,j\}}^{\infty}||S_{l-1}\partial_m{\dot{\Delta}}_j{\bar{\omega}}_n\Delta_l v_n||_{L^{\infty}}
\leq C ||{\dot{\Delta}}_j{\bar{\omega}}_n||_{L^{\infty}}\sup_{l\geq 1}||\Delta_l\nabla v_n||_{L^{\infty}}\\
&\qquad\qquad \leq C2^j ||{\dot{\Delta}}_j{\bar{v}}_n||_{L^{\infty}}\sup_{l\geq 1}||\Delta_l\nabla v_n||_{L^{\infty}}\
\end{split}
\end{equation}
where we applied Bernstein's Lemma and took the sum to get the second inequality.  We multiply ($\ref{comm2}$) by $2^{-j}$, and we take the supremum over $j\in\Z$.  This yields 
\begin{equation*}
\sup_{j\in\Z}2^{-j}||T_{\partial_m{\dot{\Delta}}_j {\bar{\omega}}_n}v_n(t)||_{L^{\infty}} \leq C||{\bar{v}}_n(t)||_{L^{\infty}}||v_n(t)||_{C^1_{\ast}}\leq C||{\bar{v}}_n(t)||_{L^{\infty}},
\end{equation*}
where we applied Remark \ref{C1revised} to get the last inequality.  We bound $||{\bar{v}}_n||_{L^{\infty}}$ using ($\ref{log}$) and ($\ref{middlefreqest}$).  We conclude that
\begin{equation*}
\sup_{j\in\Z}2^{-j}||T_{\partial_m{\dot{\Delta}}_j {\bar{\omega}}_n}v_n(t)||_{L^{\infty}} \leq C2^{-n}+Cp||{\bar{v}}_n(t)||_{{\dot{B}}^{0}_{\infty,\infty}}^{\frac{1}{q}}.
\end{equation*}
Moreover, since the Fourier transform of $S_{l-1}\partial_m{\bar{\omega}}_n\Delta_l v_n$ has support in an annulus with inner and outer radius of order $2^l$,  we have for $j\geq 0$
\begin{equation}\label{comm3}
\begin{split}
&||{\dot{\Delta}}_j(T_{\partial_m {\bar{\omega}}_n}v_n)||_{L^{\infty}} = ||{\dot{\Delta}}_j (\sum_{l\geq 1}S_{l-1}\partial_m{\bar{\omega}}_n\Delta_l v_n)||_{L^{\infty}}\\
& \leq \sum_{l=\max\{1,j-4\}}^{j+4} \sum_{k\leq l} 2^{2k}2^{-l} ||\Delta_k {\bar{v}}_n||_{L^{\infty}}||\Delta_l \nabla {v}_n||_{L^{\infty}}\leq C2^j ||{\bar{v}}_n||_{L^{\infty}}||\omega^0||_{L^{\infty}},
\end{split}
\end{equation}
where we used Bernstein's Lemma to get the first inequality, and we used Lemma \ref{CZhighfreq} and ($\ref{NSvortbound}$) to get the second inequality.  For the case $j<0$, $||{\dot{\Delta}}_j(T_{\partial_m {\bar{\omega}}_n}v_n)||_{L^{\infty}}$ is identically $0$.  Therefore ($\ref{comm3}$) still holds.  We again bound $||{\bar{v}}_n||_{L^{\infty}}$ using ($\ref{log}$) and ($\ref{middlefreqest}$), we multiply ($\ref{comm3}$) by $2^{-j}$, and we take the supremum over $j\in\Z$, which yields
\begin{equation*}
\sup_{j\in\Z}2^{-j}||{\dot{\Delta}}_j(T_{\partial_m {\bar{\omega}}_n}v_n)(t)||_{L^{\infty}} \leq C2^{-n}+Cp||{\bar{v}}_n(t)||_{{\dot{B}}^{0}_{\infty,\infty}}^{\frac{1}{q}}.
\end{equation*} 

To estimate the remainder, we again expand the commutator and consider each piece separately.  We break $v_n$ into a low-frequency term and high-frequency term, and we consider $||{\dot{\Delta}}_j(\partial_mR((Id-S_2)v_n,{\bar{\omega}}_n))||_{L^{\infty}}$.  We have 

\begin{equation*}
\begin{split}
&\sup_{j\in\Z}2^{-j}||{\dot{\Delta}}_j(\partial_mR((Id-S_2)v_n,{\bar{\omega}}_n))||_{L^{\infty}}\\
&\qquad\leq C\sum_{l\geq 0}\sum_{i=-1}^1 ||\Delta_{l-i} (Id-S_2)v_n||_{L^{\infty}}||\Delta_l {\bar{\omega}}_n||_{L^{\infty}}\\ 
&\qquad\leq C\sum_{l\geq 0}\sum_{i=-1}^12^{-(l-i)}||\Delta_{l-i}(Id-S_2)\nabla v_n||_{L^{\infty}}2^l ||\Delta_l {\bar{v}}_n||^{\frac{1}{q}}_{L^{\infty}}||\Delta_l {\bar{v}}_n||^{\frac{1}{p}}_{L^{\infty}}\\
&\qquad\leq C||{\bar{v}}_n||^{\frac{1}{q}}_{{\dot B}^{0}_{\infty,\infty}} \sum_{l\geq 0}\sum_{i=-1}^1 2^{-(l-i)}||\Delta_{l-i}(Id-S_2)\nabla v_n||_{L^{\infty}}2^l 2^{-\frac{l}{p}}||\Delta_l {\bar{\omega}}_n||^{\frac{1}{p}}_{L^{\infty}},
\end{split}
\end{equation*}
where we repeatedly used Bernstein's Lemma.  We use Lemma \ref{CZhighfreq} and ($\ref{NSvortbound}$) to bound $||\Delta_{l-i}(Id-S_2)\nabla v_n||_{L^{\infty}}$ and $||\Delta_l {\bar{\omega}}_n||_{L^{\infty}}$ by $||\omega^0||_{L^{\infty}}$ and sum over $l$ to conclude that
\begin{equation*}
\begin{split}
&\sup_{j\in\Z}2^{-j}||{\dot{\Delta}}_j(\partial_mR((Id-S_2)v_n,{\bar{\omega}}_n))||_{L^{\infty}}\\
&\qquad \qquad\leq C||{\bar{v}}_n||^{\frac{1}{q}}_{{\dot B}^{0}_{\infty,\infty}} \frac{1}{1-2^{-\frac{1}{p}}}\leq Cp||{\bar{v}}_n||^{\frac{1}{q}}_{{\dot B}^{0}_{\infty,\infty}} 
\end{split}
\end{equation*}

To bound the low frequencies, we again apply Bernstein's Lemma and the definition of the remainder term.  We write
\begin{equation*}
\begin{split}
&\sup_{j\in\Z}2^{-j}||{\dot{\Delta}}_j(\partial_mR(S_2v_n,{\bar{\omega}}_n))||_{L^{\infty}} \leq C\sum_{l\leq 3}\sum_{i=-1}^1 ||\Delta_{l-i} S_2 v_n||_{L^{\infty}}||\Delta_l {\bar{\omega}}_n||_{L^{\infty}}\\
&\qquad\qquad\qquad \leq C||v||_{L^{\infty}}||{\bar{v}}_n||_{L^{\infty}}.\\
\end{split}
\end{equation*}  
The last inequality follows from the bound $||\Delta_l{\bar{\omega}}_n||_{L^{\infty}}\leq||\Delta_l\nabla\bar{ v}||_{L^{\infty}}$, Bernstein's Lemma, and the observation that $l\leq 3$.  We now bound $||v||_{L^{\infty}}$ using Remark \ref{C1revised} and we bound $||{\bar{v}}_n||_{L^{\infty}}$ using ($\ref{log}$) and ($\ref{middlefreqest}$).  This yields
\begin{equation*}
\begin{split}
&\sup_{j\in\Z}2^{-j}||{\dot{\Delta}}_j(\partial_mR(S_2v_n,{\bar{\omega}}_n))(t)||_{L^{\infty}}\\
&\qquad\qquad\leq C2^{-n}+Cp||{\bar{v}}_n||^{\frac{1}{q}}_{{\dot B}^{0}_{\infty,\infty}}.
\end{split}
\end{equation*}
It remains to bound $\sup_{j\in\Z}2^{-j}||\partial_mR(v_n,{\dot{\Delta}}_j{\bar{\omega}}_n))||_{L^{\infty}}$.  Again we break $v_n$ into a low-frequency and high-frequency term.  We first estimate the high-frequency term.  We reintroduce the sum over $m$ and utilize the divergence-free property of $v$ to put the partial derivative $\partial_m$ on ${\bar{\omega}}_n$.  We then apply Bernstein's Lemma to conclude that for any fixed $j\in\Z$,

\begin{equation}\label{comm7}
\begin{split}
&\sum_m||R((Id-S_0)v^m_n,{\dot{\Delta}}_j\partial_m{\bar{\omega}}_n))||_{L^{\infty}}\leq C\sum_{|k-l|\leq 1} 2^{l-k}||\Delta_k\nabla v_n||_{L^{\infty}} ||{\dot{\Delta}}_j\Delta_l{\bar{\omega}}_n||_{L^{\infty}}\\
&\qquad\leq C\sup_{k\geq -1}||\Delta_k\nabla v_n||_{L^{\infty}}||{\dot{\Delta}}_j{\bar{\omega}}_n||_{L^{\infty}}\leq ||v||_{C^1_{\ast}}2^j||{\bar{v}}_n||_{L^{\infty}}.
\end{split}
\end{equation}
The second inequality above follows because for fixed $j\geq 0$, we are summing only over $l$ satisfying $|l-j|\leq 1$, while for fixed $j<0$, we are only considering $l$ satisfying $-1\leq l \leq 1$.  We now multiply ($\ref{comm7}$) by $2^{-j}$, and we take the supremum over $j\in\Z$, which yields
\begin{equation*}
\begin{split}
&\sup_{j\in\Z}2^{-j}\sum_m||\partial_mR((Id-S_0)v^m_n,{\dot{\Delta}}_j{\bar{\omega}}_n))(t)||_{L^{\infty}} \leq C||{\bar{v}}_n(t)||_{L^{\infty}}\\
&\qquad \leq C2^{-n}+Cp||{\bar{v}}_n(t)||^{\frac{1}{q}}_{{\dot B}^{0}_{\infty,\infty}},
\end{split}
\end{equation*}
where we used ($\ref{log}$) and ($\ref{middlefreqest}$) to bound $||{\bar{v}}_n(t)||_{L^{\infty}}$, and we used Remark \ref{C1revised} to bound $||v(t)||_{C^1_{\ast}}$.  For the low-frequency term, we again use the divergence-free condition on $v_n$ to write
\begin{equation*}
\begin{split}
&\sup_{j\in\Z}2^{-j}\sum_m||\partial_mR(S_0v^m_n,{\dot{\Delta}}_j{\bar{\omega}}_n)||_{L^{\infty}}\\ &\qquad\leq \sup_{j\in\Z}2^{-j}\sum_{|k-l|\leq 1}||\Delta_kS_0 v_n||_{L^{\infty}} 2^{2j}||{\dot{\Delta}}_j\Delta_l{\bar{v}}_n||_{L^{\infty}}\leq C||v||_{L^{\infty}}||{\bar{v}}_n||_{L^{\infty}}.
\end{split}
\end{equation*}
To get the first inequality, we bounded $||{\dot{\Delta}}_j\partial_m{\bar{\omega}}_n||_{L^{\infty}}$ with $||{\dot{\Delta}}_j\partial_m\nabla{\bar{v}}_n||_{L^{\infty}}$ and applied Bernstein's Lemma.  The second inequality follows from the observation that we are only considering $k\leq 1$, and therefore, by properties of our partition of unity, we are only considering $j\leq 3$.   As with previous terms, we use Remark \ref{C1revised} to bound $||v||_{L^{\infty}}$ and we use ($\ref{log}$) combined with ($\ref{middlefreqest}$) to bound $||{\bar{v}}_n||_{L^{\infty}}$.  We conclude that 
\begin{equation*}
\begin{split}
&\sup_{j\in\Z}2^{-j}\sum_m||\partial_mR(S_0v^m_n,{\dot{\Delta}}_j{\bar{\omega}}_n)(t)||_{L^{\infty}}\\
&\leq C2^{-n}+Cp||{\bar{v}}_n(t)||^{\frac{1}{q}}_{{\dot B}^{0}_{\infty,\infty}}. 
\end{split}
\end{equation*}
This completes the proof of Lemma \ref{comm}.
\end{proof}
Combining ($\ref{keyineq}$) with ($\ref{commest}$) and using the property that $A\geq 1$ gives

\begin{equation}\label{similartochemin}
\begin{split}
\partial_t \delta_n(t) 
&\leq C\left(2^{-n\alpha} +  p \delta_n(t)^{\frac{1}{q}} + \frac{1}{A}p ||{\bar{v}}_n(t)||^{\frac{1}{q}}_{{\dot B}^{0}_{\infty,\infty}}\right)\\
&\leq C\left(2^{-n\alpha} +  p \delta_n(t)^{\frac{1}{q}}\right).
\end{split}
\end{equation}

To complete the proof we closely follow an argument in \cite{Chemin}.  We let $p=2-\log \delta_n(t)$ so that ($\ref{similartochemin}$) reduces to 
\begin{equation}\label{asinchemin}
\partial_t \delta_n(t) \leq C\left(2^{-n\alpha} +  (2-\log \delta_n(t))(\delta_n(t))^{1-\frac{1}{2-\log \delta_n(t)}}\right).
\end{equation}
We now observe that 
\begin{equation*}
\delta_n(t)^{1-\frac{1}{2-\log \delta_n(t)}}\leq \delta_n(t)^{1+\frac{1}{\log \delta_n(t)}}\leq C\delta_n(t),
\end{equation*}
where we used the equality $\delta_n(t)^{\frac{1}{\log\delta_n(t)}}=e$ to obtain the last inequality.  We then see that ($\ref{asinchemin}$) reduces to 
\begin{equation*}
\partial_t \delta_n(t)\leq C\left(2^{-n\alpha} +  (2-\log \delta_n(t))\delta_n(t)\right).
\end{equation*}
Define the function $\mu$ by $\mu(r)=r(2-\log r)$.  Then, integrating with respect to time, we get 
\begin{equation}\label{almost}
\begin{split}
\delta_n(t) &\leq \delta_n(0) + C2^{-n\alpha} t +C\int_0^t{\mu(\delta_n(s))} ds\\
&\leq ||{\bar{v}}_n(0)||_{{\dot{B}}^0_{\infty,\infty}} + C2^{-n\alpha} t +C \int_0^t{\mu(\delta_n(s))} ds,
\end{split}
\end{equation}
where we used that $\delta_n(0)=\frac{||{\bar{v}}_n(0)||_{{\dot{B}}^0_{\infty,\infty}}}{A}\leq ||{\bar{v}}_n(0)||_{{\dot{B}}^0_{\infty,\infty}}$ since $A\geq 1$.  We now use the embedding $L^{\infty}\hookrightarrow {\dot{B}}^0_{\infty,\infty}$, the definition of ${\bar{v}}_n$, and Bernstein's Lemma to conclude that 
\begin{equation}\label{initial}
\begin{split}
&||{\bar{v}}_n(0)||_{{\dot{B}}^0_{\infty,\infty}} \leq C||{\bar{v}}_n(0)||_{L^{\infty}}\\
&\qquad\qquad\leq C\sum_{j\geq n}2^{-j} ||\Delta_j \omega^0||_{L^{\infty}}\leq C2^{-n}||\omega^0||_{L^{\infty}}.  
\end{split}
\end{equation}
We combine ($\ref{initial}$) with ($\ref{almost}$).  This yields
\begin{equation*}
\delta_n(t) \leq C2^{-n\alpha}( t+1) +C \int_0^t{\mu(\delta_n(s))} ds.
\end{equation*}
We now recall Osgood's Lemma.  A proof of the lemma can be found in \cite{Chemin1}.
\begin{lemma}
Let $\rho$ be a positive borelian function, and let $\gamma$ be a locally integrable positive function.  Assume that, for some strictly positive number $\beta$, the function $\rho$ satisfies
\begin{equation*}
\rho(t)\leq \beta + \int_{t_0}^t { \gamma(s)\mu(\rho(s)) } ds.
\end{equation*}
Then
\begin{equation*}
-\phi(\rho(t)) + \phi(\beta) \leq \int_{t_0}^t {\gamma(s)} ds
\end{equation*}  
where $\phi(x)=\int_x^1{\frac{1}{\mu(r)}}dr$.
\end{lemma}
We recall that we are working on a fixed time interval $[0,T]$ and we assume $n$ is large enough to ensure that $C2^{-n\alpha}( T+1) \leq e^2$.  We let $\phi$, $\rho$, $\gamma$, and $\beta$ be given by the following:\\
\\
\indent$\phi(x)=\log(2-\log x)-\log 2$, $\rho(t)=\delta_n(t)$, $\gamma(t)=C$,\\ \\
 \indent and $\beta=C2^{-n\alpha}( T+1)$.\\
\\
Applying Osgood's Lemma, we have that for any $t\leq T$,    
\begin{equation*}
-\log\left(2-\log \delta_n(t)\right) + \log \left(2-\log(C2^{-n\alpha}( T+1)\right) \leq Ct. 
\end{equation*}
Taking the exponential twice gives
\begin{equation*}
\delta_n(t) \leq e^{2-2e^{-Ct}}(C(T+1)2^{-n\alpha})^{e^{-Ct}}.
\end{equation*}
Multiplying both sides by $A$ yields the result.  This completes the proof of Lemma \ref{key}.

\end{document}